\documentclass[reqno]{amsart}

\usepackage{amsmath,amsfonts,amssymb,amsthm,hyperref,enumitem,mathdots,kbordermatrix}
\usepackage[all]{xy}

\theoremstyle{plain}
\newtheorem{lem}{Lemma}[section]
\newtheorem{cor}[lem]{Corollary}
\newtheorem{prop}[lem]{Proposition}
\newtheorem{thm}[lem]{Theorem}

\newtheorem*{mainthm}{Main Theorem}

\theoremstyle{definition}
\newtheorem{defn}[lem]{Definition}
\newtheorem{notn}[lem]{Notation}
\newtheorem{eg}[lem]{Example}
\newtheorem*{ack}{Acknowledgements}

\theoremstyle{remark}
\newtheorem{rmk}[lem]{Remark}
\newtheorem{qn}[lem]{Question}
\newtheorem{conj}[lem]{Conjecture}
\newtheorem*{introconj}{Conjecture}

\DeclareMathOperator{\spec}{Spec}
\DeclareMathOperator{\maxspec}{maxSpec}
\DeclareMathOperator{\supp}{Supp}
\DeclareMathOperator{\ann}{ann}
\DeclareMathOperator{\ap}{Ass}
\DeclareMathOperator{\depth}{depth}
\DeclareMathOperator{\height}{height}
\DeclareMathOperator{\coker}{coker}
\DeclareMathOperator{\im}{im}
\DeclareMathOperator{\cl}{Cl}
\DeclareMathOperator{\bigmod}{Mod}
\DeclareMathOperator{\fgmod}{mod}

\DeclareMathOperator{\sgn}{sgn}
\renewcommand{\geq}{\geqslant}
\renewcommand{\leq}{\leqslant}
\renewcommand{\phi}{\varphi}
\newcommand{\ds}{\displaystyle}
\newcommand{\ophom}{\operatorname{Hom}}
\newcommand{\opext}{\operatorname{Ext}}
\newcommand{\optor}{\operatorname{Tor}}
\newcommand{\Hom}[3]{\ophom_{#1}(#2,#3)}
\newcommand{\Ext}[4]{\opext_{#1}^{#2}(#3,#4)}
\newcommand{\Tor}[4]{\optor_{#1}^{#2}(#3,#4)}
\newcommand{\bc}[3]{\operatorname{\mathcal{B}}_{#1}^{#2,#3}}
\newcommand{\bbn}{\mathbb{N}}
\newcommand{\bbz}{\mathbb{Z}}
\newcommand{\sfk}{\mathsf k}
\newcommand{\p}{\mathfrak{p}}
\newcommand{\q}{\mathfrak{q}}
\newcommand{\m}{\mathfrak{m}}
\newcommand{\s}{\mathfrak{S}}
\newcommand{\ga}{\alpha}
\newcommand{\gb}{\beta}
\newcommand{\gc}{\gamma}
\newcommand{\gd}{\delta}
\newcommand{\gz}{\zeta}
\newcommand{\gl}{\lambda}
\newcommand{\vare}{\varepsilon}
\newcommand{\ol}{\overline}
\newcommand{\ul}{\underline}
\newcommand{\ti}{\widetilde}
\newcommand{\ns}{\hspace{-6pt}}

\numberwithin{equation}{lem}

\begin{document}

\bibliographystyle{amsplain}

\author{Tony Se}
\address{Department of Mathematics,
Florida A\&M University,
203 Jackson-Davis Hall,
1617 S Martin Luther King Jr.\ Blvd,
Tallahassee, FL 32307, USA}
\email{tony.se@famu.edu}

\title{Some properties of $n$-semidualizing modules}

\date{\today}


\begin{abstract} 
Let $R$ be a commutative noetherian ring. The $n$-semidualizing modules of $R$ are generalizations of
its semidualizing modules. We will prove some basic properties of $n$-semidualizing modules. Our main result
and example shows that the divisor class group of a Gorenstein determinantal ring over a field is the set of
isomorphism classes of its 1-semidualizing modules. Finally, we pose some questions about $n$-semidualizing modules.
\end{abstract}

\maketitle

\section*{Introduction}

Throughout this paper, all rings are commutative noetherian, unless stated otherwise, $\sfk$ denotes a field, and $\bbn$
is the set of nonnegative integers. Given a ring $R$, we let
$\bigmod(R)$ denote the class of all $R$-modules and $\fgmod(R)$ the class of all finitely generated $R$-modules.
We say that $C \in \fgmod(R)$ is \emph{semidualizing} if and only if $\Hom{R}{C}{C} \cong R$ and $\Ext{R}{i}{C}{C}=0$ for all $i>0$.
Semidualizing modules were first studied abstractly by Foxby \cite{f} and Golod \cite{g}, and since then by various authors.
See \cite{sw} for an introduction to the subject. In \cite[Theorem~4.2]{sw07}, Sather-Wagstaff showed that
the only semidualizing modules of a determinantal ring $R$ over $\sfk$ are $R$ and $\omega$ up to isomorphism,
where $\omega$ is the canonical module of $R$. In this paper, we consider a generalization of semidualizing modules,
called $n$-semidualizing modules. Our definition of $n$-semidualizing modules is similar, but not identical, to that of Takahashi \cite{t}.
We will show that nontrivial $n$-semidualizing modules exist for determinantal rings.

An outline of our paper is as follows. In Section~\ref{sec:gor}, we define and prove some basic properties of $n$-semidualizing modules. 
Section~\ref{sec:normal} shows that the 1-semidualizing modules of a normal domain can be found in its divisor class group.
In Section~\ref{sec:det}, we prove our main result.
\begin{mainthm}[Theorem~\ref{thm:main}]
  Let $X$ be an $n \times n$ matrix of indeterminates over $\sfk$ and $R$ the determinantal ring $\sfk[X]/(\det(X))$.
  Then the isomorphism classes in the divisor class group of $R$ are exactly those of the 1-semidualizing modules of $R$.
\end{mainthm}

Section~\ref{sec:eg} shows that the Main Theorem does not hold in general even for Gorenstein normal domains.
Finally, we indicate some open questions about $n$-semidualizing modules in Sections~\ref{sec:gor}, \ref{sec:det} and \ref{sec:eg},
in particular Conjecture~\ref{conj}.
\begin{introconj}
  Let $X$ be an $m \times n$ matrix of indeterminates over $\sfk$ and $R$ the determinantal ring
  $\sfk[X]/(I_t(X))$ with $t \leq \min(m,n)$.
  If $0 \neq [M] \in \cl(R)$, then $M$ is exactly $(m+n-2t+1)$-semidualizing. Hence $\s_0^{m+n-2t+1}(R)=\cl(R)$.
\end{introconj}

\section{Definitions and basic properties} \label{sec:gor}

Starting with the definition of $n$-semidualizing modules, we will prove some basic results about them in this section.
Most results are similar to those in \cite{sw}, but we include or sketch their proofs for completeness. In Theorem~\ref{thm:gor},
we show that if $R$ is a Gorenstein ring with $\dim(R)=d<\infty$, then any $d$-semidualizing module of $R$ is, in a sense, trivial.

\begin{defn}
  Let $R$ be a ring and $n\in \bbn$. Then $C \in \fgmod(R)$ is \emph{$n$-semidualizing} if and only if $\Hom{R}{C}{C} \cong R$
  and $\Ext{R}{i}{C}{C}=0$ for all $0<i\leq n$. We write $\s_0^n(R)$ to denote the set of isomorphism classes of $n$-semidualizing modules of $R$.
  We say that $C$ is \emph{exactly $n$-semidualizing} if and only if $[C] \in \s_0^n(R)\setminus \s_0^{n+1}(R)$.
\end{defn}

\begin{rmk} \label{rmk:homcc}
\ 
  \begin{itemize}[leftmargin=*]
    \item Let $C \in \fgmod(R)$. Then $C$ is 0-semidualizing simply when $\Hom{R}{C}{C} \cong R$.
    \item If $C \in \bigmod(R)$, then $\Hom{R}{C}{C} \cong R$ if and only if the natural map $R \to \Hom{R}{C}{C}$ is an isomorphism
    \cite[Proposition 2.2.2(a)]{sw}.
    \item Our definition of an $n$-semidualizing module differs from that in \cite[Definition~2.3]{t} in the cases $n=0,1$.
    It is this crucial difference when $n=1$ that allows us to prove Proposition~\ref{prop:cl} and Theorem~\ref{thm:main}.
  \end{itemize}
\end{rmk}

\begin{defn}
  Let $C \in \fgmod(R)$ and $m,n \in \bbn \cup \{\infty\}$. The \emph{Bass class} $\bc{C}{m}{n}(R)$
  denotes the class of all $M \in \bigmod(R)$ that satisfy the following.
  \begin{enumerate}[label=(\alph*),leftmargin=*]
    \item The evaluation map $\xi_M^C \colon C\otimes_R\Hom{R}{C}{M}\to M$ is an isomorphism.
    \item $\Ext{R}{i}{C}{M}=0$ for all $0<i\leq m$.
    \item $\Tor{i}{R}{C}{\Hom{R}{C}{M}}=0$ for all $0<i\leq n$.
  \end{enumerate}
\end{defn}

\begin{lem} \label{lem:bass}
  Let $0\to L\to M\to N\to 0$ be an exact sequence of $R$-modules. Let $m,n\in \bbn$,
  and suppose that $C \in \fgmod(R)$ has $\supp_R(C)=\spec(R)$.
  If $N \in \bc{C}{m}{n+1}$ and $M \in \bc{C}{m+1}{n}$, then $L \in \bc{C}{m+1}{n}$.
\end{lem}

\begin{proof}
  Suppose that $N \in \bc{C}{m}{n+1}$ and $M \in \bc{C}{m+1}{n}$. Applying $\Hom{R}{C}{-}$ to get the long exact sequence
  \begin{align} \label{eqn:les}
  \begin{split}
    0 &\to \Hom{R}{C}{L} \to \Hom{R}{C}{M} \to \Hom{R}{C}{N}\\
    &\to \Ext{R}{1}{C}{L} \to \Ext{R}{1}{C}{M}=0 \to \Ext{R}{1}{C}{N} \to \cdots
  \end{split}
  \end{align}
  Now applying $C \otimes_R -$ gives the following commutative diagram with exact rows.
  \[ 
    \xymatrix@C-6pt{
    C\otimes_R \Hom{R}{C}{M} \ar[r] \ar[d]_{\wr\|}^{\xi_M^C} & C\otimes_R \Hom{R}{C}{N} \ar[r] \ar[d]_{\wr\|}^{\xi_N^C}
    & C\otimes_R \Ext{R}{1}{C}{L} \ar[r] & 0\\
    M \ar[r] & N \ar[r] & 0}
 \]
  Then $C \otimes_R \Ext{R}{1}{C}{L} = 0$, so $\Ext{R}{1}{C}{L} = 0$ by \cite[Lemma~A.2.1]{sw} since $C$ has full support,
  and $\Ext{R}{i}{C}{L}=0$ for all $2\leq i\leq m+1$ by \eqref{eqn:les}.
  Since $\Tor{1}{R}{C}{\Hom{R}{C}{N}}=0$, we can complete the diagram as follows.

  \[ 
    \xymatrix@C-6pt{
    0 \ar[r] & C\otimes_R \Hom{R}{C}{L} \ar[r] \ar[d]^{\xi_L^C} & C\otimes_R \Hom{R}{C}{M} \ar[r] \ar[d]_{\wr\|}^{\xi_M^C}
    & C\otimes_R \Hom{R}{C}{N} \ar[r] \ar[d]_{\wr\|}^{\xi_N^C} & 0\\
    0 \ar[r] & L \ar[r] & M \ar[r] & N \ar[r] & 0}
 \]
  Hence $\xi_L^C$ is also an isomorphism. The long exact sequence from the first row also shows that
  $\Tor{i}{R}{C}{\Hom{R}{C}{M}}=0$ for $0<i\leq n$. Therefore, $L \in \bc{C}{m+1}{n}$.
\end{proof}

\begin{cor} \label{cor:mj}
  Let $0 \to M \to M_0 \to M_1 \to \dots \to M_j \to 0$ be an exact sequece of $R$-modules with $j \geq 1$,
  and suppose that $C \in \fgmod(R)$ has $\supp_R(C)=\spec(R)$. If $m,n\in \bbn$ and $M_i \in \bc{C}{m+j-i}{n+i}$
  for $0\leq i\leq j$, then $M \in \bc{C}{m+j}{n}$.
\end{cor}

\begin{proof}
  Break the exact sequence into short exact sequences and use Lemma~\ref{lem:bass}.
\end{proof}

\begin{lem} \label{lem:homeval}
  Let $L \in \fgmod(R)$ and $M,N \in \bigmod(R)$. Then the natural Hom evaluation map
  $\theta_{LMN} \colon L \otimes_R \Hom{R}{M}{N} \to \Hom{R}{\Hom{R}{L}{M}}{N}$ is an isomorphism if $N$ is injective.
\end{lem}

\begin{proof}
  This is \cite[Lemma~A.1.3~(2)]{sw}.
\end{proof}

\begin{cor} \label{cor:containsinj}
  Let $n\in \bbn$ and $C \in \fgmod(R)$. Then the following are equivalent.
\begin{enumerate}[label=(\roman*),leftmargin=*]
  \item \label{item:sdm} $C$ is $n$-semidualizing.
  \item \label{item:faithful} $\bc{C}{\infty}{n}$ contains a faithfully injective $R$-module.
  \item \label{item:inj} $\bc{C}{\infty}{n}$ contains every injective $R$-module.
  \item \label{item:injdim} $\bc{C}{\infty}{n-d}$ contains every $R$-module of injective dimension $\leq d$ for all $0\leq d\leq n$.
\end{enumerate}
\end{cor}

\begin{proof}
  We follow the proof in \cite[Proposition~3.1.9]{sw}.

  \ref{item:sdm} $\Rightarrow$ \ref{item:inj}: Suppose that $C$ is $n$-semidualizing and $M$ is an injective $R$-module.
  Then $\Ext{R}{i}{C}{M}=0$ for all $i>0$.

  Next, consider a free resolution $\mathcal{F} \colon \cdots \to F_2 \to F_1 \to F_0 \to 0$
  of $C$, where each $F_i$ is finitely generated. By Lemma~\ref{lem:homeval}, there is an isomorphism of complexes
  $\mathcal{F} \otimes_R \Hom{R}{C}{M} \cong \Hom{R}{\Hom{R}{\mathcal{F}}{C}}{M}$.
  Since $M$ is injective, we have $\Tor{i}{R}{C}{\Hom{R}{C}{M}} \cong \Hom{R}{\Ext{R}{i}{C}{C}}{M}$ for all $i$.
  Since $C$ is $n$-semidualizing, we have $\Tor{i}{R}{C}{\Hom{R}{C}{M}}=0$ for all $0<i\leq n$, and for $i=0$ we have
  $C \otimes_R \Hom{R}{C}{M} = \Hom{R}{R}{M} = M$. Hence $M \in \bc{C}{\infty}{n}$.

  \ref{item:faithful} $\Rightarrow$ \ref{item:sdm}: Reverse the last few arguments in \ref{item:sdm} $\Rightarrow$ \ref{item:inj}.

  \ref{item:inj} $\Rightarrow$ \ref{item:injdim}: This follows from Corollary~\ref{cor:mj} and Proposition~\ref{prop:1sdm}.

  \ref{item:injdim} $\Rightarrow$ \ref{item:inj} $\Rightarrow$ \ref{item:faithful}: Easy. See \cite[Example~A.2.3]{sw}
  for \ref{item:inj} $\Rightarrow$ \ref{item:faithful}.
\end{proof}

\begin{defn}[\protect{\cite[page~9]{sw}}]
  Let $R$ be a ring. A module $D \in \fgmod(R)$ is \emph{dualizing} if and only if it is semidualizing and has finite injective dimension.
\end{defn}

The next Theorem generalizes \cite[Corollary~4.1.9]{sw} and \cite[Lemma~5.5]{t}.

\begin{thm} \label{thm:gor}
  Let $R$ be a Gorenstein ring with $\dim(R)=d< \infty$. If $C \in \fgmod(R)$ is $n$-semidualizing with $n\geq d$,
  then $C$ is a rank 1 projective and dualizing $R$-module. In particular, if $R$ is local, then $C \cong R$.
\end{thm}

\begin{proof}
  If $C$ is $n$-semidualizing with $n\geq d$, then $C$ is $d$-semidualizing. We have $R \in \bc{C}{\infty}{0}$
  by Corollary~\ref{cor:containsinj}~\ref{item:injdim}, so the evaluation map $\xi_R^C \colon  C\otimes_R\Hom{R}{C}{R} \to R$
  is an isomorphism. Let $\m \in \maxspec(R)$. Tensoring the map $\xi_R^C$ with the residue field $\kappa(\m)$ and by counting dimension,
  we see that $C_{\m}$ is a cyclic $R_{\m}$-module. By Proposition~\ref{prop:1sdm}, $\ann_{R_{\m}}(C_{\m})=0$, so $C_{\m} \cong R_{\m}$.
  That is, $C$ is a rank 1 projective module. Hence $C$ is semidualizing by \cite[Corollary~2.2.5]{sw}, and the rest
  of the Theorem follows from \cite[Corollary~4.1.9]{sw} since $R$ is Gorenstein.
\end{proof}

Our main theorem, Theorem~\ref{thm:main}, shows that a ring $R$ with $\dim(R)=d< \infty$ may have nontrivial $n$-semidualizing modules
with $n\leq d-2$, even when $R$ is Gorenstein. So we ask the following question.

\begin{qn} \label{qn:main}
  If $R$ is a Gorenstein ring with $\dim(R)=d<\infty$, is every $(d-1)$-semidualizing module in fact semidualizing?
  Can we remove the Gorenstein assumption?
\end{qn}

\section{Normal domains} \label{sec:normal}

In this section, in anticipation of Theorem~\ref{thm:main}, we will prove Proposition~\ref{prop:cl},
which states that if $R$ is a normal domain, then the isomorphism classes of its 1-semidualizing modules are in
its divisor class group.

We will use the description of the \emph{divisor class group} $\cl(R)$ of a normal domain $R$ in \cite[pp.~261--262]{sw07}.
Let $(-)^*=\Hom{R}{-}{R}$. We say that $M \in \fgmod(R)$ is \emph{reflexive} if and only if $M \cong M^{**}$.
Then $\cl(R)$ is the set of isomorphism classes $[M]$ of reflexive $R$-modules $M$ of rank 1. As an abelian group,
the additive identity of $\cl(R)$ is $[R]$, and the group operations are given by
\[
  [M] + [N] = [(M \otimes_R N)^{**}] \quad \text{and} \quad [M]-[N]=[\Hom{R}{N}{M}].
\]

\begin{prop} \label{prop:1sdm}
  Let $C$ be a 0-semidualizing $R$-module.
  \begin{enumerate}[label=(\alph*),leftmargin=*]
    \item One has $\ann_R(C)=0$, $\supp_R(C)=\spec(R)$, $\dim_R(C)=\dim(R)$, and $\ap_R(C)=\ap_R(R)$.
    \item Given an ideal $I \subseteq R$, one has $IC=C$ if and only if $I=R$.
    \item \label{item:tf} An element $x \in R$ is $R$-regular if and only if it is $C$-regular.
  \end{enumerate}
\end{prop}

\begin{proof}
  The proof is identical to that in \cite[Proposition~2.1.16]{sw}.
\end{proof}

Part~\ref{item:seq} of the following Proposition appears in \cite[Lemma~4.8~(1)]{t}, but our proof is slightly different,
and the proof technique will resurface in the proofs of Lemma~\ref{lem:p1} and Theorem~\ref{thm:main}.

\begin{prop} \label{prop:reg}
  Let $C$ be an $(n-1)$-semidualizing $R$-module with $n \geq 1$.
  \begin{enumerate}[label=(\alph*),leftmargin=*]
    \item \label{item:seq} The sequence $x_1,\dots,x_n \in R$ is $C$-regular if and only if it is $R$-regular.
    \item \label{item:modx} If $n \geq 2$ and $x \in R$ is $R$-regular, then $C/xC$ is an $(n-2)$-semidualizing $(R/xR)$-module.
  \end{enumerate}
\end{prop}

\begin{proof}
  We follow the proof of \cite[Theorem~2.2.6]{sw} and prove part~\ref{item:modx} first. Suppose that $n \geq 2$ and
  $x \in R$ is $R$-regular. Let $\ol{R}=R/xR$ and $\ol{C}=C/xC$. By Proposition~\ref{prop:1sdm}~\ref{item:tf},
  $x$ is $C$-regular, so we have an exact sequence
  \begin{equation} \label{eqn:mult}
    0 \to C \xrightarrow{x} C \to \ol{C} \to 0.
  \end{equation}
  Applying $\Hom{R}{C}{-}$, we have $\Ext{R}{i}{C}{\ol{C}}=0$ for all $0<i<n-1$. Since $x$ is
  both $R$- and $C$-regular, we have $\Ext{\ol{R}}{i}{\ol{C}}{\ol{C}} \cong \Ext{R}{i}{C}{\ol{C}}$ for all $i \geq 0$
  by \cite[p.\ 140, Lemma 2]{m}.
  Hence $\Ext{\ol{R}}{i}{\ol{C}}{\ol{C}}=0$ for all $0<i<n-1$. The proof that $\Hom{\ol{R}}{\ol{C}}{\ol{C}} \cong \ol{R}$
  is identical to that in \cite[Theorem~2.2.6]{sw}. Therefore, $\ol{C}$ is an $(n-2)$-semidualizing $\ol{R}$-module.

  The proof of part~\ref{item:seq} is by induction. The base case is Proposition~\ref{prop:1sdm}~\ref{item:tf},
  and the induction step is given by part~\ref{item:modx}, using $x=x_1$.
\end{proof}

\begin{eg}
  Unlike \cite[Theorem~2.2.6 (c)]{sw}, if $C$ is an $n$-semidualizing $R$-module for some $n>0$ and $I$ is a proper ideal of $R$,
  we have $\depth_R(I;C) \neq \depth(I;R)$ in general. For example, let $X$ be an $m \times m$ matrix of indeterminates
  over a field $\sfk$ with $m \geq 2$, and $R=\sfk[X]/(\det(X))$. Let $\p$, respectively $\q$, be the ideal generated by the
  $(m-1)$-minors of any $m-1$ rows, respectively columns, of $X$. In Theorem~\ref{thm:main}, we will see that 
  the 1-semidualizing modules of $R$ are exactly those isomorphic to a power of $\p$ or $\q$. However,
  by \cite[Examples~(9.27)~(d)]{bv}, the only Cohen-Macaulay modules of $R$ of rank 1 are $R$, $\p$ and $\q$
  up to isomorphism.
\end{eg}

The following lemma is elementary, but we include it here for ease of reference.

\begin{lem} \label{lem:rank1}
  Let $R$ be a domain and $C \in \fgmod(R)$.
  \begin{enumerate}[label=(\alph*),leftmargin=*]
    \item If $C$ is 0-semidualizing, then it has rank 1. \label{item:rank}
    \item If $R$ is normal and $C$ has rank 1, then $C$ is 0-semidualizing. \label{item:hom}
  \end{enumerate}
\end{lem}

\begin{proof}
  Let $K$ be the quotient field of $R$. We note that an $R$-module $C$ has rank 1 if and only if
  it is isomorphic to a nonzero ideal of $R$.

  \ref{item:rank} If $\Hom{R}{C}{C} \cong R$, then tensoring with $K$
  gives $\Hom{K}{C \otimes K}{C \otimes K} \cong K$, and the result follows from counting dimension.

  \ref{item:hom} Suppose that $C\neq 0$ is isomorphic to an ideal of $R$. Then $\Hom{R}{C}{C} \subseteq K$.
  If $R$ is normal, then $\Hom{R}{C}{C} \cong R$ by the ``determinantal trick''.
\end{proof}

\begin{rmk} \label{rmk:rigid}
  Let $A$ be a ring. Recall (from algebraic geometry and representation theory) that $M \in \bigmod(A)$ is
  \emph{rigid} if and only if $\Ext{A}{1}{M}{M}=0$. Thus, by Lemma~\ref{lem:rank1} and its proof, if $R$ is a normal domain
  and $C \in \fgmod(R)$, then $[C] \in \s_0^1(R)$ if and only if $C$ is isomorphic to a nonzero rigid ideal of $R$.
\end{rmk}

\begin{prop} \label{prop:cl}
  Let $R$ be a normal domain and $C \in \fgmod(R)$. Then $[C] \in \s_0^1(R)$ if and only if $C$ is a rank 1 reflexive module
 and $\Ext{R}{1}{C}{C}=0$. In particular, $\s_0^1(R) \subseteq \cl(R)$, that is, the rigid ideals of $R$ are reflexive.
\end{prop}

\begin{proof}
  The proof is similar to that of \cite[Lemma~1.1]{SWSeSp3}. Suppose that $C$ is 1-semidualizing.
  Then $\Ext{R}{1}{C}{C}=0$ by definition, and by Lemma~\ref{lem:rank1}, $C$ has rank~1.
  To see that $C$ is reflexive, we verify the conditions in \cite[Proposition~1.4.1~(b)]{bh}.

  First, let $\p$ be a prime ideal of $R$.
  Suppose that $\height(\p)=1$. Since $R$ is $(R_1)$, the ring $R_{\p}$ is a discrete valuation ring.
  By Proposition~\ref{prop:1sdm}~\ref{item:tf}, $C_{\p}$ is torsion-free.
  By the structure theorem for principal ideal domains, $C_{\p} \cong R_{\p}$, so $C_{\p}$ is reflexive.

  Next, suppose that that $\height(\p)\geq 2$. Since $R$ is $(S_2)$, we have $\depth(R_{\p}) \geq 2$.
  Since $C$ is 1-semidualizing, we also have $\depth(C_{\p}) \geq 2$ by Proposition~\ref{prop:reg}~\ref{item:seq}.
  Therefore, $C$ is reflexive.

  Conversely, if $C$ is rank 1 reflexive, then $[C]\in \cl(R)$, so $\Hom{R}{C}{C}\cong R$.

  Finally, by Remark~\ref{rmk:rigid}, $\s_0^1(R) \subseteq \cl(R)$ if and only if the rigid ideals of $R$ are reflexive.
\end{proof}

\section{Gorenstein determinantal rings} \label{sec:det}

Our goal in this section is to prove Theorem~\ref{thm:main}, which states that the isomorphism classes of the 1-semidualizing
modules of a Gorenstein determinantal ring over a field are exactly those in the divisor class group of the ring. The Theorem also shows that
these rings give a positive answer to the first half of Question~\ref{qn:main} and Question~\ref{qn:subgp}.

Let us first review some material about determinantal rings. 
Let $\sfk$ be a field and $X=(X_{ij})$ an $m \times n$ matrix of indeterminates over $\sfk$. Let $1 < t \leq \min(m,n)$
and $I_t(X)$ be the ideal generated by all $t$-minors of $X$. Consider \emph{determinantal rings} of the form $R=R_t(X)=\sfk[X]/I_t(X)$.
Then $R$ is a Cohen-Macaulay normal domain by \cite[Remark~(2.12) and Corollary~(5.17)]{bv}, and $R$ is Gorenstein
if and only if $m=n$ by \cite[Corollary~(8.9)]{bv}. Let $\p$, respectively $\q$,
be the ideal of $R$ generated by the $(t-1)$-minors of any $t-1$ rows, respectively columns, of $X$. By \cite[Corollary~(8.4)]{bv},
$\cl(R) = \bbz[\p] = \bbz[\q]$ since $[\p]=-[\q]$, and \cite[Corollary~(7.10)]{bv} shows that $\ell[\p]=[\p^{(\ell)}]=[\p^\ell]$ and
$\ell[\q]=[\q^{(\ell)}]=[\q^\ell]$ for all $\ell \in \bbn$.

Let $[a_1,\dots,a_t \mid b_1,\dots,b_t]$ denote the determinant with rows $a_1,\dots,a_t$ and columns $b_1,\dots,b_t$ of $X$.
Let $\Pi$ be poset of $R$ consisting of the residue classes of all $t$-minors of $X$ with $t<n$, with partial order given by $[a_1,\dots,a_u \mid b_1,\dots,b_u]
\leq [c_1,\dots,c_v \mid d_1,\dots,d_v]$ if and only if $u \geq v$ and $a_1 \leq c_1$, \ldots, $a_v \leq c_v$, $b_1 \leq d_1$, \ldots, $b_v \leq d_v$
\cite[p.~46]{bv}. Then $R$ is a \emph{graded algebra with straightening law} over $\Pi$ by \cite[Theorem~(5.3)]{bv}.
The products $\gz_1\cdots \gz_\nu$ with $\nu \in \bbn$, $\gz_i \in \Pi$ and $\gz_1 \leq \cdots \leq \gz_\nu$ are called
\emph{standard monomials} \cite[p.~38]{bv}. By \cite[Proposition~(4.1)]{bv}, the standard monomials form a $\sfk$-basis of $R$.
The \emph{straightening laws} over $R$ are the relations $\gz \eta = \sum a_\mu \mu$, where $\gz,\eta \in \Pi$ are incomparable,
$0 \neq a_\mu \in \sfk$, $\mu$ is a standard monomial, and every $\mu$ has a factor $\gd \in \Pi$ such that $\gd \leq \gz$ and $\gd \leq \eta$
\cite[p.~38]{bv}.

\begin{notn}
  Let $X=(X_{ij})$ be a matrix of determinates. We then define $B_{ij}=\{X_{k\ell} \mid k=i \text{ or } \ell=j\}$, that is,
  the set of variables that are in row $i$ or column $j$.
\end{notn}

\begin{rmk} \label{rmk:review}
Let $X=(X_{ij})$ be an $m \times n$ matrix of indeterminates, $Y=X\setminus B_{mn}$ and $1<t \leq \min(m,n)$. By \cite[Proposition~(2.4)]{bv},
there is an isomorphism $R_t(X)[x_{mn}^{-1}] \cong R_{t-1}(Y)[B_{mn}][X_{mn}^{-1}]$ given by the following map.
\begin{align*}
  X_{ij} &\mapsto X_{ij} + X_{mj}X_{in}X_{mn}^{-1} \qquad \text{for all } 1\leq i \leq m-1 \text{ and } 1 \leq j \leq n-1,\\
  X_{mj} &\mapsto X_{mj},\quad X_{in} \mapsto X_{in}
\end{align*}
\end{rmk}

\begin{lem} \label{lem:localize}
  Let $X,Y,t$ be as in Remark~\ref{rmk:review}.
  Let $I(t),J(t)$ denote a power of $\p$ or $\q$ in $R_t(X)$. Let $I(t-1),J(t-1)$
  denote the corresponding powers of $\p$ or $\q$ in $R_{t-1}(Y)$. Let $S=R_{t-1}(Y)[B_{mn}][X_{mn}^{-1}]$. Then for all $i \geq 0$,
  \[ \Ext{R_t(X)}{i}{I(t)}{J(t)}_{x_{mn}} \cong \Ext{R_{t-1}(Y)}{i}{I(t-1)}{J(t-1)} \otimes_{R_{t-1}(Y)} S, \]
  and similarly for $\Tor{i}{R_t(X)}{I(t)}{J(t)}$.
\end{lem}

\begin{proof}
  First, note that the isomorphism in Remark~\ref{rmk:review} maps the ideals
  $I(t)_{x_{mn}}$ and $J(t)_{x_{mn}}$ to the extensions of $I(t-1),J(t-1)$ in $S$ respectively. We have
  \begin{align*}
    \Ext{R_t(X)}{i}{I(t)}{J(t)}_{x_{mn}}
    &\cong \Ext{R_t(X)_{x_{mn}}}{i}{I(t)_{x_{mn}}}{J(t)_{x_{mn}}}\\
    &\cong \Ext{S}{i}{I(t-1)\otimes S}{J(t-1)\otimes S}\\
    &\cong \Ext{R_{t-1}(Y)}{i}{I(t-1)}{J(t-1)}\otimes_{R_{t-1}(Y)} S,
  \end{align*}
  where the last two isomorphisms hold since $S$ is faithfully flat over $R_{t-1}(Y)$.
\end{proof}

\begin{lem}[\protect{\cite[Lemma~4.4]{bv}}] \label{lem:plue}
  Consider an $m \times p$ matrix over a commutative ring with $m \leq p$ and indices
  $c_1,\dots,c_k,e_{\ell},\dots,e_m,d_1,\dots,d_s \in \{1,\dots,p\}$ such that $s=2m-k-(m-\ell+1)>m$
  and $u=m-k>0$. Then we have
  \[
    \sum_{\substack{i_1<\cdots<i_u\\ i_{u+1}<\cdots<i_s\\ \{1,\dots,s\}=\{i_1,\dots,i_s\}}}
    \sgn(i_1,\dots,i_s)[c_1,\dots,c_k,d_{i_1},\dots,d_{i_u}][d_{i_{u+1}},\dots,d_{i_s},e_\ell,\dots,e_m]=0.
  \]
\end{lem}

\begin{notn} \label{notn:tix}
  For the rest of this section, we let $X=(X_{ij})$ be an $n \times n$ matrix of determinates over $\sfk$, $R=R_n(X)$,
  $M_{ij}$ the $(i,j)$-minor of $X$, $C_{ij}$ the $(i,j)$-cofactor of $X$, and $x_{ij},m_{ij},c_{ij}$
  the images of $X_{ij},M_{ij},C_{ij}$ in $R$ respectively.
  As in \cite[pp.~45--46]{bv}, we let $\ti{X}$ be an $n\times 2n$ matrix by adding $n$ columns of indeterminates
  to the right of $X$, and consider the epimorphism $\sfk[\ti{X}] \to \sfk[X]$ given by mapping the entries in
  $\ti{X}$ to the corresponding entry in the matrix
  \[
    \left( \begin{matrix}
      X_{11} & \cdots & X_{1n} & 0 & \cdots & \cdots & 0 & 1\\
      &&&  \vdots && \iddots & \iddots & 0\\
      \vdots && \vdots & \vdots & \iddots & \iddots & \iddots & \vdots\\
      &&& 0 & \iddots & \iddots && \vdots\\
      X_{n1} & \cdots & X_{nn} & 1 & 0 & \cdots & \cdots & 0
    \end{matrix} \right).
  \]
\end{notn}

\begin{cor}
  Let $j_0 \in \{1,\dots,n-1\}$, $t \in \{j_0,\dots,n-1\}$, $1\leq a_1 < \cdots < a_t \leq n$,
  $j_0<b_{j_0+1}< \cdots < b_t \leq n$. Then in $R=R_n(X)$ we have
  \begin{equation} \label{eqn:cofactor}
    \sum_{1\leq j \leq n} c_{nj} [a_1,\dots,a_t \mid j,1,2,\dots,j_0-1,b_{j_0+1},\dots,b_t]=0.
  \end{equation}
\end{cor}

\begin{proof}
  Apply Lemma~\ref{lem:plue} to the matrix $\ti{X}$ over $\sfk[\ti{X}]$ with $m=n$, $p=2n$,
  $k=0$, $\ell=2$, $s=n+1$, $u=n$, $d_1=1$, \ldots, $d_{n+1}=n+1$, $e_2=1$, \ldots, $e_{j_0}=j_0-1$,
  $e_{j_0+1}=b_{j_0+1}$, \ldots, $e_n=b_n$, where
  $\{a_1,\dots,a_t,2n+1-b_n,\dots,2n+1-b_{t+1}\}=\{1,\dots,n\}$, to get
  \[
    \sum_{1\leq j \leq n+1} (-1)^{n+1-j} [1,\dots,j-1,j+1,\dots,n,n+1][j,1,2,\dots,j_0-1,b_{j_0+1},\dots,b_n]=0.
  \]
  Apply the epimorphism $\sfk[\ti{X}] \to \sfk[X]$ in Notation~\ref{notn:tix} and then the natural map
  $\sfk[X] \to R$ to get
  \[
    \sum_{1\leq j \leq n} (-1)^{n+1-j} m_{nj}[a_1,\dots,a_t \mid j,1,2,\dots,j_0-1,b_{j_0+1},\dots,b_t]=0,
  \]
  and note that $c_{nj}=(-1)^{n+j}m_{nj}$.
\end{proof}

\begin{prop}[\protect{\cite[Example 4.1]{hw}}] \label{prop:periodic}
  Let $n\geq 2$. Consider the matrices
  \begin{gather*}
    \ti{\ga} = \left( \begin{matrix}
        X_{11} & -X_{21} & \cdots & (-1)^n X_{n-1,1} & (-1)^{n+1}X_{n1}\\
        -X_{12} & X_{22} & \cdots & (-1)^{n+1} X_{n-1,2} & (-1)^{n+2}X_{n2}\\
        \vdots & & \vdots & & \vdots\\
        (-1)^n X_{1,n-1} & (-1)^{1+n} X_{2,n-1} & \cdots & X_{n-1,n-1} & -X_{n,n-1}\\
        (-1)^{1+n}X_{1n} & (-1)^{2+n} X_{2n} & \cdots & -X_{n-1,n} & X_{nn}
    \end{matrix} \right)\\
    \text{and}\quad \ti{\gb} = \left( \begin{matrix}
        M_{11} & M_{12} & \cdots & M_{1,n-1} & M_{1n}\\
        M_{21} & M_{22} & \cdots & M_{2,n-1} & M_{2n}\\
        \vdots & & \vdots & & \vdots\\
        M_{n1} & M_{n2} & \cdots & M_{n,n-1} & M_{nn}
    \end{matrix} \right)
  \end{gather*}
  over $\sfk[X]$. Let $\p = (m_{n1},m_{n2},\dots,m_{nn})$.
  Let $\ga=((-1)^{j+i}x_{ji})$ and $\gb=(m_{ij})$ over $R$. Then the complex
  \begin{equation} \label{eqn:resolution}
      \cdots \xrightarrow{\gb} R^{\oplus n} \xrightarrow{\ga} R^{\oplus n}
      \xrightarrow{\gb} R^{\oplus n} \xrightarrow{\ga} R^{\oplus n}
  \end{equation}
  of period 2 is a free resolution of $\p$.
\end{prop}

\begin{proof}
  First, $(\ti{\ga},\ti{\gb})$ is a matrix factorization of $\det(X)$, so \eqref{eqn:resolution} is a free resolution of
  $\coker \ti{\ga} = \coker \ga$ \cite[Proposition 5.1]{e}. Let us augment \eqref{eqn:resolution} with $R^{\oplus n} \xrightarrow{\vare} \p \to 0$,
  where $\vare$ is given by the matrix $( \begin{matrix} m_{n1} & m_{n2} & \cdots & m_{n,n-1} & m_{nn} \end{matrix} )$.
  We need to show that $\ker \vare = \im \ga$. Certainly $\im \ga \subseteq \ker \vare$ by expanding $\det(X)$.

  We need to show that $\im \ga$ generates $\ker \vare$. Let $\Psi = \{ m_{n1},\dots,m_{nn}\} \subset \Pi$,
  so that $\Psi$ is an ideal of $\Pi$, i.e.\ if $\gz \in \Psi$ and $\eta \leq \gz$, then $\eta \in \Psi$ \cite[p.~50]{bv}.
  Since $m_{nj}$ is the residue class of $[1,\dots,n-1 \mid 1,\dots,j-1,j+1,\dots,n]$, we have $m_{n1}>m_{n2}>
  \cdots > m_{nn}$.
  Let $e_j$ be the basis element of $R^{\oplus n}$ such that $\vare(e_j)=m_{nj}$. By \cite[Proposition~(5.6) (b)]{bv},
  we have $\ker \vare = \im \ga$ once we show that $\im \ga$ contains elements
  \[
    g_{\xi j} = \xi e_j - \sum_{j<k\leq n} r_{\xi k j} e_k \text{ with } r_{\xi k j} \in R
  \]
  for all $\xi \in \Pi$ and $j \in \{1,\dots,n-1\}$ such that $\xi \ngeqslant m_{nj}$.

  Let $\xi$ be the residue class of
  $[a_1,\dots,a_t \mid b_1,\dots,b_t]$, where $t \in \{1,\dots,n-1\}$, $a_1<\cdots<a_t$ and $b_1<\cdots<b_t$.
  If $\xi \ngeqslant m_{nj}$, then we have $t \geq j$ and $b_1=1,\dots,b_j=j$. If $t=1$, then we simply let
  the $g_{x_{i1} 1}$ be given by the columns of $\ga$, that is,
  \[
    g_{x_{i1} 1} = \sum_{1\leq k \leq n} (-1)^{i+k} x_{ik} e_k.
  \]
  If $2\leq t \leq n-1$, then \eqref{eqn:cofactor} gives
  \[
    \sum_{1\leq k \leq n} (-1)^{n+k} m_{nk}[a_1,\dots,a_t \mid k,1,2,\dots,j-1,b_{j+1},\dots,b_t]=0
  \]
  in $R$. Of course, the first $j-1$ terms are simply 0. Let
  \[
    g_{\xi j} = \sum_{j\leq k \leq n} (-1)^{n+k}[a_1,\dots,a_t \mid k,1,2,\dots,j-1,b_{j+1},\dots,b_t] e_k,
  \]
  so that $g_{\xi j} \in \ker \vare$. To see that $g_{\xi j} \in \im \ga$, expand the $t$-minors along the first column to get
  \begin{align*}
    g_{\xi j} &= \sum_{1\leq k \leq n} (-1)^{n+k} \sum_{1\leq i \leq t} (-1)^{i+1} x_{a_i k}\\
    &\qquad \qquad [a_1,\dots,a_{i-1},a_{i+1},\dots,a_t \mid 1,2,\dots,j-1,b_{j+1},\dots,b_t]e_k\\
    &= \sum_{1\leq i \leq t} (-1)^{n+i+1-a_i} [a_1,\dots,a_{i-1},a_{i+1},\dots,a_t \mid 1,2,\dots,j-1,b_{j+1},\dots,b_t]\\
    &\qquad \qquad \sum_{1\leq k\leq n} (-1)^{a_i+k} x_{a_i k}e_k.
  \end{align*}
  Hence $g_{\xi j} \in \im \ga$ by the case $t=1$, and the proof is complete.
\end{proof}

\begin{lem} \label{lem:resolution2}
  Let $n=2$, $\p=(x_{12},x_{11})=(m_{21},m_{22})$ and $\ell \geq 1$. Then the complex
  \begin{equation} \label{eqn:pell}
      \cdots \xrightarrow{\ga^{\oplus \ell}} R^{\oplus 2\ell} \xrightarrow{\gb^{\oplus \ell}} R^{\oplus 2\ell} \xrightarrow{\ga^{\oplus \ell}} R^{\oplus 2\ell}
      \xrightarrow{\gb^{\oplus \ell}} R^{\oplus 2\ell} \xrightarrow{\gc} R^{\oplus (\ell + 1)}
  \end{equation}
  is a free resolution of $\p^\ell$ with period 2 after the map $\gc$ when $\ell >1$, where
  \[
    \ga = \left( \begin{matrix}
      x_{11} & -x_{21}\\
      -x_{12} & x_{22}
    \end{matrix} \right), \qquad
    \gb = \left( \begin{matrix}
      x_{22} & x_{21}\\
      x_{12} & x_{11}
    \end{matrix} \right),
  \]
  $\ga^{\oplus \ell},\gb^{\oplus \ell}$ are the matrix direct sums of $\ga,\gb$ respectively,
  and $\gc$ is given by the $(\ell + 1)\times 2\ell$ matrix with a copy of $\ga$ starting from entries
  $(1,1),(2,3),\dots,(\ell, 2\ell-1)$ and 0s in all other entries.
\end{lem}

\begin{proof}
  If $\ell=1$, then \eqref{eqn:pell} is simply \eqref{eqn:resolution}. So let us consider the case $\ell > 1$.

  First, we need to show that $R^{\oplus 2\ell} \xrightarrow{\gc} R^{\oplus(\ell + 1)} \xrightarrow{\vare} \p^{\ell}$
  is a presentation of $\p^{\ell}$, where $\vare$ is the natural projection map onto
  $\p^\ell = (x_{12}^\ell, x_{12}^{\ell-1}x_{11},\dots,x_{12}x_{11}^{\ell-1},x_{11}^\ell)$, and $\gc$ is described
  as in the Lemma. For example, when $\ell=3$, we have
  \[
    \gc = \left( \begin{matrix}
      x_{11} & -x_{21} & 0 & 0 & 0 & 0\\
      -x_{12} & x_{22} & x_{11} & -x_{21} & 0 & 0\\
      0 & 0 & -x_{12} & x_{22} & x_{11} & -x_{21}\\
      0 & 0 & 0 & 0 & -x_{12} & x_{22}
    \end{matrix} \right).
  \]
  Certainly $\im \gc \subseteq \ker \vare$. Conversely, suppose that $\ul{r}=(r_0,\dots,r_\ell)^T \in \ker \vare$.
  Reducing $\ul{r}$ modulo $\im \gc$, we may assume that the terms in $r_0$ involve $x_{12},x_{22}$ only.
  Since $\vare(\ul{r})=0$, we have $r_0 x_{12}^\ell \in x_{11}R$, so $r_0=0$ since there are no more relations
  in $R$ to obtain a factor of $x_{11}$ from $r_0 x_{12}$. Similarly, $r_1,\dots,r_{\ell -1}=0$. Finally, if
  $r_\ell x_{11}^\ell=0$, then $r_\ell=0$ since $R$ is a domain. Therefore, $\ul{r} \in \im \gc$, and hence $\ker \vare = \im \gc$.

  It remains to show that the sequence $R^{\oplus 2\ell} \xrightarrow{\gb^{\oplus \ell}} R^{\oplus 2\ell} \xrightarrow{\gc} R^{\oplus (\ell + 1)}$
  is exact, since the rest of \eqref{eqn:pell} is exact by the exactness of \eqref{eqn:resolution}.
  Certainly $\im \left( \gb^{\oplus \ell} \right) \subseteq \ker \gc$. Conversely, suppose that $\ul{r}=(r_1,\dots,r_{2\ell})^T \in \ker \gc$.
  From row 1 of $\gc$ we get $x_{11}r_1 - x_{21}r_2=0$. As in \eqref{eqn:resolution}, the sequence
  \[
      \cdots \xrightarrow{\ga} R^{\oplus 2} \xrightarrow{\gb} R^{\oplus 2}
      \xrightarrow{\ga} R^{\oplus 2} \xrightarrow{\gb} R^{\oplus 2}
  \]
  is a free resolution of $\q=(x_{11},x_{21})$ when augmented by $\vare'=( \begin{matrix} x_{11} & -x_{21} \end{matrix} )$.
  Then $(r_1,r_2)^T\in \ker \vare' = \im \gb=\ker \ga$, so $-x_{12}r_1+x_{22}r_2=0$ as well.
  Since $\ul{r} \in \ker \gc$, from row 2 of $\gc$ we get $(r_3,r_4)^T \in \ker \vare' = \im \gb$, and so on.
  Therefore, $\ul{r} \in \im \left( \gb^{\oplus \ell} \right)$, and hence $\ker \gc = \im \left( \gb^{\oplus \ell} \right)$.
\end{proof}

\begin{prop} \label{prop:2by2}
   Let $n=2$, $\p=(x_{12},x_{11})=(m_{21},m_{22})$ and $\ell \geq 1$. Then $\Ext{R}{i}{\p^\ell}{\p^\ell} \neq 0$
  if and only if $i$ is even, where $i \geq 0$. In particular, the ideal $\p^\ell$ is exactly 1-semidualizing, and $\s_0^1(R)=\cl(R)$.
\end{prop}

\begin{proof}
  When $i=0$, we have $\Hom{R}{\p^\ell}{\p^\ell}\cong R$ since $[\p^\ell]=\ell [\p] \in \cl(R)$.

  For $i \neq 0$, apply $\Hom{R}{-}{\p^\ell}$ to \eqref{eqn:pell} to get
  \begin{gather*}
    (\p^\ell)^{\oplus (\ell + 1)} \xrightarrow{\gc^T} (\p^\ell)^{\oplus 2\ell} \xrightarrow{(\gb^T)^{\oplus \ell}} (\p^\ell)^{\oplus 2\ell}
    \xrightarrow{(\ga^T)^{\oplus \ell}} (\p^\ell)^{\oplus 2\ell} \xrightarrow{(\gb^T)^{\oplus \ell}}
    (\p^\ell)^{\oplus 2\ell} \xrightarrow{(\ga^T)^{\oplus \ell}} \cdots,\\
    \ga^T = \left( \begin{matrix}
      x_{11} & -x_{12}\\
      -x_{21} & x_{22}
    \end{matrix} \right) \quad \text{ and } \quad
    \gb^T = \left( \begin{matrix}
      x_{22} & x_{12}\\
      x_{21} & x_{11}
    \end{matrix} \right).
  \end{gather*}
  Then $(x_{12}^\ell,x_{11}x_{12}^{\ell - 1},\dots,x_{12}^\ell,x_{11}x_{12}^{\ell - 1})^T \in \ker ((\ga^T)^{\oplus \ell}|(\p^\ell)^{\oplus 2\ell})
  \setminus \im((\gb^T)^{\oplus \ell}|(\p^\ell)^{\oplus 2\ell})$. Hence $\Ext{R}{i}{\p^\ell}{\p^\ell}\neq 0$ if $i$ is even.

  Next, we show that $\ker((\gb^T)^{\oplus \ell} |(\p^\ell)^{\oplus 2\ell}) \subseteq \im ((\ga^T)^{\oplus \ell} | (\p^\ell)^{\oplus 2 \ell})$.
  Suppose that $\ul{r'} \in \ker((\gb^T)^{\oplus \ell} |(\p^\ell)^{\oplus 2\ell})$.
  Since $(\ti{\ga},\ti{\gb})$ in Proposition~\ref{prop:periodic} is a matrix factorization, we have $\ker \gb^T = \im \ga^T$ (in $R^{\oplus 2}$).
  So $\ul{r'} \in \im (\ga^T)^{\oplus \ell} \cap (\p^\ell)^{\oplus 2\ell}$. Let $\ul{r'} = (\ga^T)^{\oplus \ell}(\ul{r})$ with
  $\ul{r}=(r_1,r_2,\dots,r_{2\ell-1},r_{2\ell})^T$. We need to show that $\ul{r} \in (\p^\ell)^{\oplus 2\ell}$.
  Since $(\ga^T)^{\oplus \ell}(\ul{r}) \in (\p^\ell)^{\oplus 2\ell}$, we have $-x_{2,1}r_{2j-1}+x_{2,2}r_{2j} \in \p^\ell$ for $j=1,\dots,\ell$.
  Reducing $\ul{r}$ modulo $\ker (\ga^T)^{\oplus \ell}$ and noting that $\ker \ga^T = \im \gb^T$,
  we may assume that the terms in $x_{2,2}r_{2j}$ involve $x_{1,2},x_{2,2}$ only. Since these terms do not appear in $x_{2,1}r_{2j-1}$,
  we have $r_{2j} \in x_{1,2}^\ell R$. Then $x_{2,1}r_{2j-1} \in \p^\ell$. Now for each term $\mu$ in $r_{2j-1}$,
  the total degree of $x_{1,1}$ and $x_{1,2}$ in $x_{2,1}\mu$ is well-defined in $R$.
  So $r_{2j-1} \in \p^\ell$ and hence $\ul{r} \in (\p^\ell)^{\oplus 2\ell}$.
  Therefore, if $i$ is odd and $i \neq 1$, then $\Ext{R}{i}{\p^\ell}{\p^\ell}= 0$.

  It remains to show that $\ker ((\gb^T)^{\oplus \ell} |(\p^\ell)^{\oplus 2\ell}) \subseteq \im (\gc^T|(\p^\ell)^{\oplus (\ell + 1)})$ (when $\ell > 1$),
  or $(\ga^T)^{\oplus \ell}((\p^\ell)^{\oplus 2\ell}) \subseteq \gc^T((\p^\ell)^{\oplus (\ell + 1)})$. For example, when $\ell=3$, we have
  \[
    \gc^T = \left( \begin{matrix}
      x_{11} & -x_{12} & 0 & 0\\
      -x_{21} & x_{22} & 0 & 0\\
      0 & x_{11} & -x_{12} & 0\\
      0 & -x_{21} & x_{22} & 0\\
      0 & 0 & x_{11} & -x_{12}\\
      0 & 0 & -x_{21} & x_{22}
    \end{matrix} \right).
  \]
  Let $\{e_j \mid j = 1,\dots,2\ell\}$ be the standard basis of $R^{\oplus 2\ell}$. Consider the columns
  \[
    g_{1j}=x_{11}e_{2j-1}-x_{21}e_{2j} \text{ and } g_{2j}=-x_{12}e_{2j-1}+x_{22}e_{2j}
  \]
  of $(\ga^T)^{\oplus \ell}$, where $j=1,\dots,\ell$, and
  \[
    h_1=g_{11},\ h_{\ell + 1}=g_{2\ell}, \text{ and } h_j=-x_{12}e_{2j-3}+x_{22}e_{2j-2}+x_{11}e_{2j-1}-x_{21}e_{2j}
  \]
  of $\gc^T$, where $j=2,\dots,\ell$. We need to show that $x_{11}^\nu x_{12}^{\ell - \nu} g_{kj} \in \gc^T((\p^\ell)^{\oplus (\ell + 1)})$
  for all $k=1,2$, $j=1,\dots,\ell$ and $\nu=0,1,\dots,\ell$. This is true because
  \begin{itemize}[leftmargin=*]
    \item $g_{11}=h_1$ and $g_{2\ell} = h_{\ell + 1}$,
    \item $x_{11}^\nu x_{12}^{\ell - \nu} g_{1j_0} = \ds \sum_{1 \leq j \leq j_0} x_{11}^{\nu-j_0+j} x_{12}^{\ell - \nu+j_0-j} h_j$
    for all $2\leq j_0 \leq \ell$ and $j_0-1 \leq \nu \leq \ell$,
    \item $x_{11}^\nu x_{12}^{\ell - \nu} g_{1j_0} = \ds \sum_{j_0 < j \leq \ell+1} -x_{11}^{\nu-j_0+j} x_{12}^{\ell - \nu+j_0-j} h_j$
    for all $2\leq j_0 \leq \ell$ and $0\leq \nu < j_0$,
  \end{itemize}
  and similarly for $x_{11}^\nu x_{12}^{\ell - \nu} g_{2j_0}$, and we have shown that $\Ext{R}{1}{\p^\ell}{\p^\ell}= 0$.

  By symmetry, we have $\Ext{R}{i}{\q^\ell}{\q^\ell}\neq 0$ if and only if $i$ is even for all $\ell \geq 1$ and $i \geq 0$,
  and hence $\s_0^1(R)=\cl(R)$ by Proposition~\ref{prop:cl}.
\end{proof}

\begin{lem} \label{lem:mnn}
  Let $n \geq 2$, $\ell \geq 0$, $\p=(m_{n1},\dots,m_{nn})$ and $\ol{\p^\ell}=\p^\ell/x_{nn}\p^\ell$, with
  $\p^0=R$. Then:
  \begin{enumerate}[label=(\alph*),leftmargin=*]
    \item \label{item:nzd} $m_{nn}$ is a nonzero divisor on $\ol{\p^\ell}$.
    \item \label{item:colon} $(x_{nn}\p^\ell :_{\p^{\ell-k}}m_{nn}^k)=x_{nn}\p^{\ell-k}$ for all $0\leq k \leq \ell$. In particular,
    $\ann_R(\ol{\p^\ell}) = x_{nn}R$.
    \item \label{item:exchange} Let $\ell \geq 2$, $\nu \in \{1,\dots,\ell-1\}$ and $r,s \in \p^k$ with $k \geq \ell - \nu$.
    If $\ol{m_{nn}^\nu s} = \ol{m_{n,n-1}m_{nn}^{\nu-1} r}$,
    then $r=r_1+m_{nn}r_2$ for some $r_1 \in x_{nn}\p^k$ and $r_2 \in \p^{k-1}$ such that $s-m_{n,n-1}r_2 \in x_{nn}\p^k$.
  \end{enumerate}
\end{lem}

\begin{proof}
  \ref{item:nzd} Let $r,s \in \p^\ell$ be such that $m_{nn}r = x_{nn}s$. Write $r,s$ as a linear combination of standard monomials over $\Pi$.
  Let $\Psi=\{m_{n1},\dots,m_{nn}\}$ as in the proof of Proposition~\ref{prop:periodic}.
  Since $\Psi$ is an ideal of $\Pi$, each standard monomial in $r,s$ is in $\p^\ell$ by the argument of \cite[Proposition~4.1]{bv}. Since $m_{nn},x_{nn}$
  are the smallest and largest elements in $\Pi$ respectively, no straightening laws are used when writing $m_{nn}r,x_{nn}s$ in terms of
  standard monomials. Then each standard monomial in $r$ is in $x_{nn}\p^\ell$, so $r \in x_{nn}\p^\ell$. Part~\ref{item:colon} is similar.

  \ref{item:exchange} Write $r$ as a linear combination of standard monomials over $\Psi$, and let $r'$ consist of the terms that have a factor
  of $m_{nn}$, and $r_1$ be the rest of the terms, so that $r=r_1+m_{nn}r_2$, $r_1 \in \p^k$ and $r_2 \in \p^{k-1}$, where $r'=m_{nn}r_2$.
  Then we have
  \[
    m_{nn}^\nu s - m_{n,n-1}m_{nn}^{\nu-1} r = m_{nn}^\nu (s - m_{n,n-1}r_2) - m_{nn}^{\nu-1} m_{n,n-1} r_1 \in x_{nn}\p^\ell.
  \]
  Since $m_{n,n-1}$ is the smallest element in $\Psi \setminus \{m_{nn}\}$, no straightening laws are used when writing
  $m_{nn}^{\nu-1} m_{n,n-1} r_1$ in terms of standard monomials, and no standard monomial in $m_{nn}^{\nu-1} m_{n,n-1} r_1$
  is a multiple of $m_{nn}^\nu$. Therefore, $r_1 \in x_{nn}\p^k$ and $s-m_{n,n-1}r_2 \in x_{nn}\p^k$.
\end{proof}

\begin{lem} \label{lem:p1}
  Let $n>2$ and $\p=(m_{n1},\dots,m_{nn})$. Then $\Ext{R}{i}{\p}{\p} \neq 0$ if and only if $i$ is even, where $i \geq 0$.
  In particular, the ideal $\p$ is exactly 1-semidualizing.
\end{lem}

\begin{proof}
  Apply $\Hom{R}{-}{\p}$ to the resolution \eqref{eqn:resolution} to get
  \[
    \p^{\oplus n} \xrightarrow{\ga^T} \p^{\oplus n} \xrightarrow{\gb^T} \p^{\oplus n}
    \xrightarrow{\ga^T} \p^{\oplus n} \xrightarrow{\gb^T} \p^{\oplus n} \to \cdots
  \]
  Then $(m_{n1},m_{n2},\dots,m_{nn})^T \in \ker (\ga^T|\p^{\oplus n}) \setminus \im(\gb^T|\p^{\oplus n})$, so
  $\Ext{R}{i}{\p}{\p}\neq 0$ if $i$ is even.

  Since \eqref{eqn:resolution} has period 2, it remains to show that $\Ext{R}{1}{\p}{\p}=0$.
  Consider the exact sequence \eqref{eqn:mult} with $C=\p$ and $x=x_{nn}$. Apply $\Hom{R}{\p}{-}$, using the
  notation in Proposition~\ref{prop:reg}, to get the commutative diagram
  \begin{equation} \label{eqn:olp}
  \begin{split}
    \xymatrix{ 0 \ar[r] & R \ar[r]^{x_{nn}} \ar[d]_{\wr\|} & R \ar[r] \ar[d]_{\wr\|} & \ol{R} \ar[r] \ar[d] & 0\\
    0 \ar[r] & \Hom{R}{\p}{\p} \ar[r]^{x_{nn}} & \Hom{R}{\p}{\p} \ar[r] & \Hom{R}{\p}{\ol{\p}}  }
  \end{split}
  \end{equation}
  where the vertical maps are the natural maps.
  Let us show that $\ol{R} \cong \Hom{R}{\p}{\ol{\p}}$. Apply $\Hom{R}{-}{\ol{\p}}$ to \eqref{eqn:resolution} to get
  \[
    \ol{\p}^{\oplus n} \xrightarrow{\ga^T} \ol{\p}^{\oplus n} \to \cdots,
  \]
  so that $\Hom{R}{\p}{\ol{\p}} \cong \ker(\ga^T|\ol{\p}^{\oplus n})$. We need to solve the system of equations
  $\ga^T(\ol{r_1},\dots,\ol{r_n})^T = (0,\dots,0)^T$ with $r_1,\dots,r_n \in \p$, that is,

  \begin{alignat*}{10}
    x_{11}\ol{r_1} & {}- x_{12}\ol{r_2} &{}+ \cdots &{}+  (-1)^{1+n}x_{1n}\ol{r_n} &&=  0\\
    -x_{21}\ol{r_1} &{}+ x_{22}\ol{r_2} &{}+ \cdots &{}+ (-1)^{2+n}x_{2n}\ol{r_n} &&=  0\\
    &&&&&\hspace{4pt}\vdots\\
    (-1)^{n+1}x_{n1}\ol{r_1} &{}+  (-1)^{n+2}x_{n2}\ol{r_2} &{}+ \cdots &{}+ x_{nn}\ol{r_n} &&=  0
  \end{alignat*}
  Let $\rho_k$ denote the $k$th equation. Then for any $j_0=1,\dots,n-1$,
  \[
    \sum_{1\leq k \leq n-1} [1,\dots,k-1,k+1,\dots,n-1 \mid 1,\dots,j_0-1,j_0+1,\dots,n-1] \rho_k
  \]
  gives $\ol{m_{nn}r_{j_0}\vphantom{P}}-\ol{m_{nj_0}r_n\vphantom{P}}=0$. In particular,
  $\ol{m_{nn}r_{n,n-1}\vphantom{P}}=\ol{m_{n,n-1}r_n\vphantom{P}}$. Then by Lemma~\ref{lem:mnn}~\ref{item:exchange},
  we get $\ol{r_n\vphantom{P}}=r'_n \ol{m_{nn}\vphantom{P}}$ and $\ol{r_{n,n-1}\vphantom{P}}=r'_n \ol{m_{n,n-1}\vphantom{P}}$
  for some $r'_n \in R$. In general, $m_{nn}\ol{r_{j_0}} - m_{nn} \ol{r'_n m_{nj_0}} = 0$,
  and so $\ol{r_{j_0}\vphantom{P}} = r'_n \ol{m_{nj_0}\vphantom{P}}$ by Lemma~\ref{lem:mnn}~\ref{item:nzd}.
  Hence $(\ol{r_1},\dots,\ol{r_n})=r'_n(\ol{m_{n1}\vphantom{P}},\dots,\ol{m_{nn}\vphantom{P}})$,
  and $\ker(\ga^T|\ol{\p}^{\oplus n}) = R(\ol{m_{n1}\vphantom{P}},\dots,\ol{m_{nn}\vphantom{P}})^T \cong \ol{R}$
  by Lemma~\ref{lem:mnn}~\ref{item:colon}.

  We now have $\Hom{\ol{R}}{\ol{\p}}{\ol{\p}} \cong \Hom{R}{\p}{\ol{\p}} \cong \ol{R}$. By Remark~\ref{rmk:homcc},
  the vertical map on the right in \eqref{eqn:olp} is an isomorphism, so the bottom row is exact. Continuing the
  long exact sequence shows that the map $\Ext{R}{1}{\p}{\p} \xrightarrow{x_{nn}} \Ext{R}{1}{\p}{\p}$
  is injective. By Lemma~\ref{lem:localize}, Proposition~\ref{prop:2by2} and induction on $n \geq 2$,
  we have $\Ext{R}{1}{\p}{\p}_{x_{nn}} = 0$, that is, $x_{nn} \in \sqrt{\ann_R(\Ext{R}{1}{\p}{\p})}$.
  Therefore, $\Ext{R}{1}{\p}{\p}=0$.
\end{proof}

\begin{defn} \label{defn:gamma}
  Let $n\geq 2$ and $\ell \geq 1$. We define a $\binom{n+\ell-1}{\ell} \times n\binom{n+\ell-2}{\ell-1}$
  matrix $\gc_\ell$ as follows. The rows of $\gc_\ell$ are labeled by the standard monomial generators $\gl$ of $\p^\ell$ in lexicographic order
  and the columns by ordered pairs $(\mu,j)$, where $\mu$ is a standard monomial generator of $\p^{\ell - 1}$ and $j \in \{1,\dots, n\}$,
  first in lexicographic order in $\mu$, then in ascending order in $j$. For each $\mu$,
  we place a copy of $\ga$ as in Proposition~\ref{prop:periodic} at the minor at columns $(\mu,1),\dots,(\mu,n)$ and
  rows $\mu m_{n1},\dots,\mu m_{nn}$. The rest of the entries of $\gc_\ell$ are 0. For example, if $n=3$, then $\gc_2$ is the matrix
  \renewcommand{\kbldelim}{(}
  \renewcommand{\kbrdelim}{)}
  \setlength{\kbcolsep}{0pt}
  \[
    \kbordermatrix{
    & \ns (m_{31},1) \ns & \ns (m_{31},2) \ns & \ns (m_{31},3) \ns & \ns (m_{32},1) \ns & \ns (m_{32},2) \ns & \ns (m_{32},3) \ns
    & \ns (m_{33},1) \ns & \ns (m_{33},2) \ns & \ns (m_{33},3) \ns \\
    m_{31}^2 & x_{11} & -x_{21} & x_{31} & 0 & 0 & 0 & 0 & 0 & 0 \\
    m_{31}m_{32} & -x_{12} & x_{22} & -x_{32} & x_{11} & -x_{21} & x_{31} & 0 & 0 & 0 \\
    m_{31}m_{33} & x_{13} & -x_{23} & x_{33} & 0 & 0 & 0 & x_{11} & -x_{21} & x_{31} \\
    m_{32}^2 & 0 & 0 & 0 & -x_{12} & x_{22} & -x_{32} & 0 & 0 & 0 \\
    m_{32}m_{33} & 0 & 0 & 0 & x_{13} & -x_{23} & x_{33} & -x_{12} & x_{22} & -x_{32} \\
    m_{33}^2 & 0 & 0 & 0 & 0 & 0 & 0 & x_{13} & -x_{23} & x_{33}
    }.
  \]
  When $n=2$, the definition of these matrices agrees with that of $\gc$ in Lemma~\ref{lem:resolution2}, and when $\ell=1$,
  we have $\gc_1=\ga$.

  The $\gc_\ell$ can also be defined inductively. Let $\gc_1=\ga$.
  Assuming that $\gc_{\ell - 1}$ has been defined, for each $m_{nj}$ with $j \in \{1,\dots,n\}$, place a copy of $\gc_{\ell - 1}$
  at the minor at the columns $(\mu,1),\dots,(\mu,n)$ and rows $\gl$ of $\gc_\ell$ where $m_{nj}$ is a factor of $\mu$ and $\gl$.
  Then put 0 in the rest of the entries of $\gc_\ell$.
\end{defn}

\begin{thm} \label{thm:main}
   Let $\p=(m_{n1},\dots,m_{nn})$ and $\ell \geq 1$. Then the ideal $\p^\ell$ is exactly 1-semidualizing, and $\s_0^1(R)=\cl(R)$.
\end{thm}

\begin{proof}
  The case for $n=2$ is in Proposition~\ref{prop:2by2} and the case for $\ell = 1$ in Lemma~\ref{lem:p1}.
  For $n>2$ and $\ell > 1$, we will first show that the matrix $\gc_\ell$ in Notation~\ref{defn:gamma} gives
  a finite presentation
  \begin{equation} \label{eqn:pres}
    R^{\oplus n \binom{n+\ell-2}{\ell - 1}} \xrightarrow{\gc_\ell} R^{\oplus \binom{n+\ell-1}{\ell}} \xrightarrow{\vare} \p^\ell
  \end{equation}
  of $\p^\ell$. Order the standard monomial generators of $\p^\ell$ in lexicographic order, and let $\vare$ be the
  natural projection map. Certainly $\im \gc_\ell \subseteq \ker \vare$. Conversely, let $\gl$ range over the standard monomial
  generators of $\p^\ell$, let $\ul{r}=(r_\gl)^T \in \ker \vare$, and write each $r_\gl$ as a linear combination of standard monomials.
  The elements $g_{\xi j}$ in the proof of Proposition~\ref{prop:periodic} show that the straightening relations
  that involve a factor $m_{nj}$ of $\gl$ are generated by columns $(\gl/m_{nj},1),\dots,(\gl/m_{nj},n)$ of $\gc_\ell$.
  Then modulo $\im \gc_\ell$, we may assume that no straightening relations are used when finding $\vare(\ul{r})=\sum_\gl r_\gl \gl$.
  Thus, for each $\gl$, the factors of the standard monomials that appear in $r_\gl$ are all $\geq$ those in $\gl$.
  Since $\vare(\ul{r})=0$, we have $\ul{r}=\ul{0}$ modulo $\im \gc_\ell$. Therefore, $\im \gc_\ell = \ker \vare$.

  Following the proof of Lemma~\ref{lem:p1}, we apply $\Hom{R}{-}{\ol{\p^\ell}}$ to \eqref{eqn:pres} truncated at $\vare$,
  so that $\Hom{R}{\p^\ell}{\ol{\p^\ell}} \cong \ker \left( \left. \gc_\ell^T\right| \ol{\p^\ell}^{\oplus \binom{n+\ell-1}{\ell}} \right)$,
  and show that the latter is isomorphic to $\ol{R}=R/x_{nn}R$. The proof is by induction on $1\leq \nu \leq \ell$ that
  \[
    \ker \left( \left. \gc_\nu^T\right| \ol{\p^\ell}^{\oplus \binom{n+\nu-1}{\nu}} \right) =
    \left\{ \sum_\eta \ol{r \eta \vphantom{P}}e_\eta \mid r \in \p^{\ell - \nu} \right\},\]
  where $\eta$ runs through the standard monomial generators of $\p^\nu$, and $\{e_\eta\}$ is the standard basis
  of $R^{\oplus \binom{n+\nu-1}{\nu}}$. When $\nu=1$,
  the proof of Lemma~\ref{lem:p1} shows that $\ker(\ga^T|\ol{\p^\ell}^{\oplus n}) =
  \{ (\ol{r m_{n1}\vphantom{h}},\dots,\ol{r m_{nn}\vphantom{h}})^T \mid r \in \p^{\ell -1}\}$. Let $1<\nu \leq \ell$,
  and let $\ul{s}=(\ol{s_\eta}) \in \ol{\p^\ell}^{\oplus \binom{n+\nu-1}{\nu}}$ be such that $\gc_\nu^T(\ul{s})=\ul{0}$.
  Let $\mu$ run through the standard monomial generators of $\p^{\nu-1}$. For each $j \in \{1,\dots,n\}$,
  apply the induction hypothesis to rows $(\mu,1),\dots,(\mu,n)$ and columns $\eta$ of $\gc_\ell^T$,
  where $m_{nj}$ is a factor of $\mu$ and $\eta$, to get
  \[
    \sum_{m_{nj} \mid \eta} \ol{s_{\eta}\vphantom{P}}e_{\eta}
    = \sum_{m_{nj} \mid \eta} \ol{r_j (\eta/m_{nj})\vphantom{P}}e_{\eta}
  \]
  for some $r_j \in \p^{\ell-\nu+1}$. Then $j=n-1,n$ gives us
  \[
    \ol{s_{m_{n,n-1}m_{nn}^{\nu-1}}\vphantom{P}}
    = \ol{r_{n-1}m_{nn}^{\nu-1}}
    = \ol{r_n m_{n,n-1}m_{nn}^{\nu-2}}.
  \]
  By Lemma~\ref{lem:mnn}~\ref{item:exchange}, there are $t_{n-1} \in x_{nn} \p^{\ell - \nu + 1}$ and
  $t_n \in \p^{\ell - \nu}$ such that $r_n=t_{n-1}+m_{nn}t_n$ and $r_{n-1}-m_{n,n-1}t_n \in x_{nn}\p^{\ell - \nu + 1}$.
  Then $\ol{s_\eta \vphantom{t}} = \ol{t_n \eta}$ whenever $m_{nn} \mid \eta$ or $m_{n,n-1} \mid \eta$.
  For $j \neq n-1,n$ we have
  \[
    \ol{s_{m_{nj}m_{nn}^{\nu-1}}\vphantom{P}}
    = \ol{r_j m_{nn}^{\nu-1}}
    = \ol{r_n m_{nj}m_{nn}^{\nu-2}}
    = \ol{t_n m_{nj}m_{nn}^{\nu - 1}}.
  \]
  Lemma~\ref{lem:mnn}~\ref{item:colon} shows that $r_j- m_{nj} t_n \in x_{nn} \p^{\ell - \nu + 1}$, so
  $\ol{s_\eta \vphantom{t}} = \ol{t_n \eta}$ whenever $m_{nj} \mid \eta$. The induction is now complete,
  and the case $\nu = \ell$ and Lemma~\ref{lem:mnn}~\ref{item:colon} show that
  $\Hom{R}{\p^\ell}{\ol{\p^\ell}} \cong \ol{R}$.

  The rest of the argument in Lemma~\ref{lem:p1}, using \eqref{eqn:olp} with $\p^\ell$ instead of $\p$,
  shows that $\Ext{R}{1}{\p^\ell}{\p^\ell} = 0$. Lemma~\ref{lem:localize} and Proposition~\ref{prop:2by2}
  with induction on $n$ show that $\Ext{R}{2}{\p^\ell}{\p^\ell} \neq 0$. Hence $\p^\ell$ is exactly 1-semidualizing.
  By symmetry, the ideals $\q^\ell$ are also exactly 1-semidualizing, and hence $\s_0^1(R)=\cl(R)$.
\end{proof}

\begin{cor}
  Let $n \geq 2$, $R=R_n(X)$ and $d=\dim R$. Then any $(d-1)$-semidualizing module of $R$ is semidualizing.
  The result is sharp with $n=2$.
\end{cor}

\begin{proof}
  This follows from Theorem~\ref{thm:main}, since $d=n^2-1$.
\end{proof}

\begin{eg}
  When $n>2$, we do not necessarily have $\Ext{R}{i}{\p^\ell}{\p^\ell}=0$ for all odd $i$. For example,
  let $n=3$ and $\m$ be the homogeneous maximal ideal of $R$. Then $\dim R=8$ and $\depth_\m \p^3 = 6$
  by \cite[Examples~(9.27)~(d)]{bv}. Hence a minimal resolution of $\p^3$ becomes periodic of period 2 after 2 steps;
  see \cite[Theorem~6.1]{e}. A calculation with Macaulay2 \cite{M2} shows that $\Ext{R}{i}{\p^3}{\p^3}=0$
  for $i=1$ only. Lemma~\ref{lem:localize} and Theorem~\ref{thm:main} then show that for all $n \geq 3$,
  $\Ext{R}{i}{\p^3}{\p^3}=0$ for $i=1$ only.
\end{eg}

\begin{rmk}
  By Remark~\ref{rmk:rigid} and Proposition~\ref{prop:cl},
  Theorem~\ref{thm:main} states that the rigid ideals of $R$ are exactly the reflexive ideals of $R$.
  See, however, Conjecture~\ref{conj}.
\end{rmk}

\begin{conj} \label{conj}
  Let $X$ be an $m \times n$ matrix of indeterminates over $\sfk$ and $R=R_t(X)$ with $t \leq \min(m,n)$.
  If $0 \neq [M] \in \cl(R)$, then $M$ is exactly $(m+n-2t+1)$-semidualizing. Hence $\s_0^{m+n-2t+1}(R)=\cl(R)$.
  In particular, if $t=2$ and $d=\dim R = m+n-1$, then $\s_0^{d-2}(R)=\cl(R)$.
\end{conj}

\section{Another example} \label{sec:eg}

In Section~\ref{sec:normal}, we saw that $\s_0^1(R) \subseteq \cl(R)$.
Now in contrast to Theorem~\ref{thm:main}, we will show that $\s_0^1(R) \neq \cl(R)$ in general even for
Gorenstein normal domains.

\begin{eg} \cite[p.~168]{m}
  Let $\sfk$ be a field of characteristic 0, $n \geq 0$, and $R=\sfk[X,Y,Z]/(XY-Z^n)$. Then $\cl(R) \cong \bbz/n\bbz$
  with generator $[\p]$, where $\p=(x,z)$. One can check that $\p^{(m)}=(x,z^m)$ for all $0<m<n$, and that
  $\p^{(m)}$ has a free resolution
  \[
      \cdots \xrightarrow{\gb} R^{\oplus 2} \xrightarrow{\ga} R^{\oplus 2}
      \xrightarrow{\gb} R^{\oplus 2} \xrightarrow{\ga} R^{\oplus 2}
  \]
  of period 2, where
  \[
    \ti{\ga} = \left[ \begin{matrix} Y & Z^m\\ -Z^{n-m} & -X \end{matrix} \right]
    \quad \text{and} \quad
    \ti{\gb} = \left[ \begin{matrix} X & Z^m\\ -Z^{n-m} & -Y \end{matrix} \right].
  \]
  Here $(\ti{\ga},\ti{\gb})$ is a matrix factorization of $XY-Z^n$. Apply $\Hom{R}{-}{\p^{(m)}}$ to get
  \begin{gather*}
    \big(\p^{(m)} \big)^{\oplus 2} \xrightarrow{\ga^T} \big(\p^{(m)} \big)^{\oplus 2} \xrightarrow{\gb^T}
    \big(\p^{(m)} \big)^{\oplus 2} \xrightarrow{\ga^T} \big(\p^{(m)} \big)^{\oplus 2} \xrightarrow{\gb^T} \cdots
    \quad \text{where}\\
    \ga^T = \left[ \begin{matrix} y & -z^{n-m}\\ z^m & -x \end{matrix} \right]
    \quad \text{and} \quad
    \gb^T = \left[ \begin{matrix} x & -z^{n-m}\\ z^m & -y \end{matrix} \right]
  \end{gather*}
  and note that the sequence given by $\Hom{R}{-}{R}$ is exact, since $(\ti{\ga},\ti{\gb})$ is a
  matrix factorization. We have $\Hom{R}{\p^{(m)}}{\p^{(m)}} \cong R$ since $[\p^{(m)}]=m[\p] \in \cl(R)$, and
  we observe the following.
  \begin{itemize}[leftmargin=*]
    \item If $0<m \leq n/2$ and $i>0$, then $\Ext{R}{i}{\p^{(m)}}{\p^{(m)}}$ has generator
    $(z^{n-m},x)^T$ when $i$ is odd, and $\Ext{R}{i}{\p^{(m)}}{\p^{(m)}}$ has generator
    $(x,z^m)^T$ when $i$ is even.
    \item If $n/2 <m <n$ and $i>0$, then $\Ext{R}{i}{\p^{(m)}}{\p^{(m)}}$ has generator $(z^m,xz^{2m-n})^T$ when $i$ is odd,
    and $\Ext{R}{i}{\p^{(m)}}{\p^{(m)}}$ has generator $(x,z^m)^T$ when $i$ is even.
  \end{itemize}

  We see that $\s_0^1(R)=\{[R]\}$, and $\s_0^1(R) = \cl(R)$ only when $n=1$.
\end{eg}

\begin{qn} \label{qn:subgp}
  If $R$ is a noetherian normal domain, is $\s_0^1(R)$ a subgroup of $\cl(R)$?
\end{qn}

\begin{ack}
  The author thanks Olgur Celikbas and Arash Sadeghi for comments on earlier versions of the manuscript.
  Ryo Takahashi first defined $n$-semidualizing modules in \cite{t}. The author thanks him for providing feedback
  which enriched this paper.
\end{ack}


\begin{thebibliography}{10}
\bibitem{bh}
W.~Bruns and J.~Herzog, \emph{Cohen-Macaulay Rings}, Cambridge
Studies in Advanced Mathematics, \textbf{39}. Cambridge University Press,
Cambridge, 1993.

\bibitem{bv}
W.\ Bruns and U.\ Vetter, \emph{Determinantal rings}, Lecture Notes in
Mathematics, vol.\ 1327, Springer-Verlag, Berlin, 1988.

\bibitem{e}
D.\ Eisenbud, \emph{Homological algebra on a complete intersection, with an application to group representations},
Trans.\ Am.\ Math.\ Soc.\ \textbf{260} (1980) 35--64.

\bibitem{f}
H.-B.\ Foxby, \emph{Gorenstein modules and related modules}, Math.\ Scand. \textbf{31} (1972), 267--284 (1973).

\bibitem{g}
E.\ S.\ Golod, \emph{$G$-dimension and generalized perfect ideals}, Trudy Mat.\ Inst.\ Steklov. \textbf{165} (1984), 62--66,
Algebraic geometry and its applications.

\bibitem{M2} D.\ R.\ Grayson and M.\ E.\ Stillman,
\textit{Macaulay2}, a software system for research in algebraic geometry,
available at \url{http://www.math.uiuc.edu/Macaulay2/}.

\bibitem{hw} C.\ Huneke, R.\ Wiegand, \emph{Tensor products of modules and the rigidity of Tor},
Math.\ Ann. \textbf{299} (1994) 449--476; Correction: Math.\ Ann.\ \textbf{388} (2007) 291--293.

\bibitem{m}
H.\ Matsumura, \emph{Commutative ring theory}, second ed., Studies in Advanced Mathematics,
vol.\ 8, University Press, Cambridge, 1989.

\bibitem{sw}
S.~Sather-Wagstaff, \emph{Semidualizing modules}, course notes.
\url{https://ssather.people.clemson.edu/DOCS/sdm.pdf}

\bibitem{sw07}
S.~Sather-Wagstaff, \emph{Semidualizing modules and the divisor class group},
Illinois J.\ Math.\ \textbf{51} (2007), no.~1, 255--285.

\bibitem{SWSeSp3}
S.~Sather-Wagstaff, T.~Se, S.~Spiroff, \textit{Generic constructions and semidualizing modules},
Algebr.\ Represent.\ Theory \textbf{24} (2021), no.\ 4, 1071--1081. DOI: 10.1007/s10468-020-09979-5

\bibitem{t}
R.~Takahashi, \textit{A new approximation theory which unifies spherical and
Cohen–Macaulay approximations}, J.\ Pure Appl.\ Algebra \textbf{208} (2007), no.\ 2, 617--634.
\end{thebibliography}
\end{document}